\documentclass[10pt]{amsart}

\usepackage[all]{xy}
\def\pf{\begin{proof}}
\def\epf{\end{proof}}

\usepackage[mathscr]{eucal}
\usepackage[active]{srcltx}   
\usepackage{amsmath,amsthm,amssymb,amsxtra,amscd, amsbsy,eucal,color}

\newcommand{\ee}{{ [e]}}

\newcommand{\Id}{{\rm{Id}}}
\newcommand{\I}{{\mathcal I}}

\newcommand{\B}{{\mathcal B}}

\newcommand{\de}{\delta}
\newcommand{\ot}{\otimes}
\newcommand{\af}{\alpha}

\newcommand{\p}{{\mathscr{G}_0}}

\newcommand{\A}{\mathcal{A}}
\newcommand{\G}{\mathscr{G}}
\newcommand{\h}{\mathscr{H}}
\newcommand{\m}{{}^{-1}}
\newcommand{\K}{\Bbbk}

\newcommand{\vone}{\vspace{0.09cm}}

\newtheorem{teo1}{Theorem}[section]

\newtheorem{lem1}[teo1]{Lemma}
\newtheorem{cor1}[teo1]{Corollary}
\newtheorem{prop1}[teo1]{Proposition}
\theoremstyle{remark}
\newtheorem{exe1}[teo1]{Example}
\theoremstyle{remark}
\newtheorem{remark}[teo1]{Remark}
\theoremstyle{remark}
\newtheorem{notation}[teo1]{Notation}
\def\pf{\begin{proof}}
\def\epf{\end{proof}}

\begin{document}

\title[Partial actions of groupoids]{On the separability of the partial skew groupoid ring}
\author[Bagio, Pinedo]{Dirceu Bagio and Hector Pinedo}

\address{ Departamento de Matem\'atica, Universidade Federal de Santa Maria, 97105-900\\
Santa Maria-RS, Brasil}\email{bagio@smail.ufsm.br}

\address{Escuela de Matem\'aticas, Universidad Industrial de Santander, Cra. 27 Calle 9 UIS Edificio 45\\  Bucaramanga, Colombia}\email{hpinedot@uis.edu.co}

\thanks{
D. Bagio was partially supported by CAPES.
H. Pinedo was supported by FAPESP}

\dedicatory{To Antonio Paques, on his 70th birthday.}


\begin{abstract} Given a partial (resp. a global) action $\af $ of a connected finite groupoid $\G$ on a ring $\A,$ we determine necessary and sufficient conditions for the partial (resp. global) skew groupoid ring $\A\star_\af \G$ to be a separable extension of $\A$. 
\end{abstract}
\maketitle

\section*{Introduction}
We recall that an extension of associative rings $R\subseteq S$ is  {\em separable} if the multiplication map $m_S:S\otimes_R S\to R$ is a $(R,R)$-bimodule epimorphism which splits. Separable extensions of rings have been extensively studied in the literature; (see e.~g. \cite{CAS,CIM,DEM,HS,I,L1}) and more recently in \cite{NOP}.

One interesting homology property of ring extensions is semisimplicity. We recall that $R\subset S$ is semisimple if any exact sequence in the category $_S\mathcal{M}$ of left $S$-modules that splits in $_R\mathcal{M}$ also splits in $_S\mathcal{M}$. It is known that separable extensions are semisimple extensions; see e.~g.  \cite[ Remark 1.1]{CIM}.

Given a finite group $G$ such that the order of $G$ is invertible on a ring $\A$, it is well-known that the group ring $\A[G]$ is a separable extension of $\A$.  If $\G$ is a groupoid, then the groupoid ring $\A[\G]$ is a separable extension of $\A$ if and only if the orders of the associated isotropy groups are invertible elements in $\A$; see \cite{L2}. 
If the group $G$ acts non trivially on $\A$ then we have associated the skew group ring $\A\star G$. Several authors investigated questions related to the separability of $\A \subset \A\star G$; see e.~g. \cite{AS, AS1, I}. For a (twisted) partial action $\af$ of  $G$ on $\A$, necessary and sufficient conditions to the separability of the (twisted) partial skew group ring $\A\star_{\af} G$ over $\A$ were given in \cite{BLP}.

Let $\G$ be a finite connected groupoid acting partially on a ring $\A$ and  $\A\star_{\af}\G$ the corresponding partial skew groupoid ring. In this article, we give
necessary and sufficient conditions to the separability of $\A\subset \A\star_{\af}\G$; see Theorem \ref{sepa}. We also determine in $\S\,$4 under which conditions the extension $\A\subset \A\star_{\af}\G$ is separable for a global action $\alpha$ of $\G$ on $\A$.

The paper is organized as follows. After the introduction, in Section 1 we recall basic facts on groupoids and partial actions of groupoids. Moreover, we introduce trace maps that are fundamental to our purposes. Some properties of the partial skew groupoid ring are studied in Section 2. In Section 3, given a partial action $\af$ of a groupoid $\G$ on a ring $A$, we investigate when $\A\star_{\af}\G$ is a separable extension of $\A$. Firstly, in Proposition \ref{prop:redu}, we reduce this problem to connected groupoids. After that, we present our main result in Theorem \ref{sepa}. The last Section explores the same problem for global actions of groupoids.

\subsection*{Conventions}\label{subsec:conv}
In this work, {\it ring} means a  non-necessarily unital associative ring. For a subset $\I$ of an ring $\A$ we write $\I\unlhd \A$ to  denote that $\I$ is a two-sided ideal of 
 $\A$. The center of a ring $\A$ will be denoted by $C(\A)$.

\section{Partial actions of groupoids}
 
 \subsection{Groupoids}
A {\it groupoid} is a small category where every morphism is an isomorphism. We denote a groupoid by 
$s,t:\G\rightrightarrows \p$, or simply by $\G$, where $\G$ denotes the set of morphisms, $\p$ the set of objects and $s,t$ are {\it source} and {\it target} maps, respectively. We recall
that the product $gh$ exists if and only if $s(g)=t(h)$, $g,h\in \G$. In this case, $t(gh)=t(g)$ and 
$s(gh)=s(h)$. The set $\G^2$ consists  of pairs $(g,h)\in \G\times \G$ such that $s(g)=t(h)$.
Given $e,f\in\p$, we denote by $\G(e,f)$ the set of arrows from $e$ to $f$ and  by $\G(e)=\G(e,e)$ the isotropy group. Note that $\G$ induces the following equivalence relation on $\p$:
\begin{equation}\label{rel} e\sim f\,\,\, \text{if and only if}\,\,\,\G(e,f)\neq\emptyset.\end{equation} The groupoid $\G$ is {\it connected} if 
$e\sim f$, for all $e,f\in \p$. Let  $[e]$ be the
equivalence class of $e\in \p$. Throughout this paper, we assume that $\G$ is a groupoid such that $\p/\hspace*{-0.15cm}\sim$ is finite. Also,  we fix a transversal $[e_1],\ldots, [e_n]$ of 
the equivalence relation $\sim$, that is, $\p/\hspace*{-0.15cm}\sim=\{[e_1],\ldots,[e_n]\}$.
\vone
 
Given $[e]\in \p/\hspace*{-0.15cm}\sim$, consider the full subgroupoid $\G_{[e]}\rightrightarrows [e]$ of $\G$, that is, $\G_{[e]}(f,f')=\G(f,f')$, for all $f,\,f'\in [e]$. Observe that
$\G_{[e]}=\G_{[e_j]}$ whenever $e\sim e_j$. By construction, $\G_{[e]}$ is connected and it is called {\it the connected component of $\G$ associated to $e$}. Note that
\[\G=\bigsqcup\limits_{\ee\in \p/\sim}\G_{[e]}\,\,=\,\,\bigsqcup\limits_{i=1}^n\G_{[e_i]},\]
that is, $\G$ is a disjoint union of its connected components.

 \subsection{Partial actions of groupoids} Let $\G$ be a groupoid and $\A$ a ring. We recall from \cite{BP} that a {\it partial action} of $\G$ on $\A$ is a collection
 $\af=(\A_g,\af_g)_{g\in \G}$, where for each $g\in \G$, $\A_g \unlhd \A_{t(g)},$  $\A_{t(g)}\unlhd \A$, $\af_g\colon \A_{g\m}\to \A_g$ is a ring isomorphism, and the following conditions hold:
 \begin{enumerate}\renewcommand{\theenumi}{\roman{enumi}}   \renewcommand{\labelenumi}{(\theenumi)}
 	\item $\af_e$ is the identity map $\Id_{\A_e}$ of $\A_e$, for all $e\in\p$;
 	\item $\af^{-1}_h(\A_{g\m}\cap \A_h)\subseteq \A_{(gh)\m}$, for all $(g,h)\in \G^2$;
 	\item $\af_g(\af_h(x))=\af_{gh}(x)$, for all $(g,h)\in \G^2$ and $x\in \af^{-1}_h(\A_{g\m}\cap \A_h)$.
 \end{enumerate}
 Conditions \rm{(ii)} and \rm{(iii)} mean that $\af_{gh}$ is an extension of $\af_g\af_h,$ for all $(g,h)\in \G^2$. We say that $\af$ is a {\it global action} if 
 $\af_g\af_h=\af_{gh}$,  for all $(g,h)\in \G^2$. By  \cite[Remark 1.2]{BP}, $\alpha_{(e)}=(\A_g,\alpha_g)_{g\in \G(e)}$ is a partial
 action of the  group $\G(e)$ on the subring $\A_e$ of $\A$, for each $e\in\p$.
 \vone
 
 
 A partial action $\alpha=(\A_g,\alpha_g)_{g\in \G}$ of a groupoid $\G$ on a ring $\A$ is called {\it unital} if
 $\A_g=\A1_g$, where $1_g$ is a central idempotent of $\A$, for all
 $g\in \G$.\smallbreak
 
 The next result was proved in \cite[Lemma 1.1]{BP} and
 we recall it here for
 the reader's convenience.
 
 \begin{lem1}\label{BP1} Let $\af=(\A_g,\af_g)_{g\in \G}$ be a partial action of a groupoid $\G$ on a ring $\A$.
 	Then:
 	\begin{enumerate}\renewcommand{\theenumi}{\roman{enumi}}   \renewcommand{\labelenumi}{(\theenumi)}
 		\item $\af$ is global if and only if $\A_g=\A_{t(g)}$, for all $g\in \G$;
 		\item $\af_{g\m}=\af\m_g$, for all $g\in \G$;
 		\item $\af_g(\A_{g\m}\cap \A_h)=\A_{g}\cap \A_{gh}$, for all $(g,h)\in \G^2$. \qed
 	\end{enumerate}
 \end{lem1}
 \begin{remark}
 	Let $\af=(\A_g,\af_g)_{g\in \G}$ be a unital partial action of a groupoid $\G$ on  a ring $\A$,  $(g,h)\in \G^2$ and $a\in\A$. 
 	Then 
 	\begin{align}\label{iguali}
 	&\alpha_g(\alpha_h(a1_{h^{-1}})1_{g^{-1}})=\alpha_{gh}(a1_{(gh)^{-1}})1_g.&
 	\end{align} 
 	In fact, by Lemma \ref{BP1} (iii), $1_h1_{g^{-1}}=\alpha_h(1_{h^{-1}}1_{(gh)^{-1}})$ and $1_{gh}1_g=\alpha_{gh}(1_{h^{-1}}1_{(gh)^{-1}})$, for all $(g,h)\in \G^2$. Thus
 	\begin{align*}
 	\alpha_g(\alpha_h(a1_{h^{-1}})1_{g^{-1}})&=\alpha_g(\alpha_h(a1_{h^{-1}})\alpha_h(1_{h^{-1}}1_{(gh)^{-1}}))\\
 	&=\alpha_g(\alpha_h(a1_{h^{-1}}1_{(gh)^{-1}}))=\alpha_{gh}(a1_{h^{-1}}1_{(gh)^{-1}})\\
 	&=\alpha_{gh}(a1_{(gh)^{-1}})\alpha_{gh}(1_{h^{-1}}1_{(gh)^{-1}})=\alpha_{gh}(a1_{(gh)^{-1}})1_g.\
 	\end{align*}
 	
 \end{remark}

 Let $\G$ be a groupoid, $\h$ a subgroupoid of $\G$ and $\A$ a ring. Suppose that $\af=(\A_g,\af_g)_{g\in \G}$ is a partial action of $\G$ on $\A$. We denote by $\af|_{\h}$ the partial action of $\h$ on $\A$ obtained by restriction of $\alpha$ to $\h$, i.e. $\af|_{\h}=(\A_h,\af_h)_{h\in \h}$. 
 
 \begin{prop1}\label{partial_action_comp_connected}
 	Let $\G$ be a groupoid and $\A$ a ring. Then $\G$ acts partially on $\A$ if and only if $\G_{[e]}$ acts partially on $\A$, for all $[e]\in \p/\hspace*{-0.15cm}\sim$.
 \end{prop1}
 \pf
 Suppose that $\af^{(i)}=(\A^{(i)}_{g},\af^{(i)}_{g})_{g\in \G_{[e_i]}}$ is a partial action of $\G_{[e_i]}$ on $\A$, for all $i=1,\ldots,n$. For each $g\in \G$, we consider the unique $j\in\{1,\ldots,n\}$ such that $g\in \G_{[e_j]}$. We define $\A_g:=\A^{(j)}_{g}$ and $\af_g:=\af^{(j)}_{g}$. It is clear that $\af=(\A_g,\af_g)_{g\in \G}$ is a partial action of $\G$ on $\A$.
 Conversely,  if $\af$ is a partial action of $\G$ on $\A$ then  $\af|_{\G_{[e_i]}}$ is a partial action of $\G_{[e_i]}$ on $\A$, for all $i=1,\ldots, n$.\epf
 
 Let $\G$ be a groupoid and $\A$ a ring such that $\A=\oplus_{e\in \p}\A_e$, where $\A_e$ is an ideal of $\A$, for all $e\in \p$. Note that, for each $e\in \p$, we have that $\A_\ee:=\oplus_{f\in \ee}\A_f$ is an ideal of $\A$ and 
 $\A=\oplus_{i=1}^n\,\A_{[e_i]}$.

 \begin{prop1}\label{partial_action_comp_connected_2}
 	Let $\G$ be a groupoid and $\A$ a ring.
 	\begin{enumerate}\renewcommand{\theenumi}{\roman{enumi}}   \renewcommand{\labelenumi}{(\theenumi)}
 		\item Let $\af=(\A_g,\af_g)_{g\in \G}$ be a partial action of  $\G$ on $\A$ and assume that  $\A=\oplus_{e\in \p}\A_e$. Then $\alpha|_{\G_{[e]}}$ is a partial action of $\G_{[e]}$ on $\A_{[e]}$, for all $[e]\in \p/\hspace*{-0.15cm}\sim$.
 		\item Suppose that $\A=\oplus_{e\in \p}\A_e$ where $\A_e$ is an ideal of $\A$, for all $e\in\p$. If $\G_{[e]}$ acts partially on $\A_{[e]}$ for all $[e]\in \p/\hspace*{-0.15cm}\sim$, then $\G$ acts partially on $\A$.
 	\end{enumerate}
 \end{prop1}
 \pf (i) Let $[e]\in \p/\hspace*{-0.15cm}\sim$ and $g\in \G_{[e]}$. Since $\A_{s(g)}$ and $\A_{t(g)}$ are ideals of $\A_{[e]}$, it is immediate to check that $\alpha|_{\G_{[e]}}$ is a partial action of $\G_{[e]}$ on $\A_{[e]}$. 
 \smallbreak
 \noindent (ii) Assume that $\af^{(i)}=(\A^{(i)}_g,\af^{(i)}_g)_{g\in \G_{[e_i]}}$ is a partial action of $\G_{[e_i]}$ on $\A_{[e_i]}$, for all $i=1,\ldots,n$. Given $g\in \G$, there exists a unique $j\in\{1,\ldots,n\}$ such that $g\in \G_{[e_j]}$. Consider $\A_g:=\A^{(j)}_g$, $\af_g:=\af^{(j)}_g$  and $\af=(\A_g,\af_g)_{g\in\G}$. It is easy to verify that $\af$ is a partial action of $\G$ on $\A$.\epf

 \begin{notation} Let $\G$, $\A$ and $\alpha$ be as in Proposition \ref{partial_action_comp_connected_2} (i). We shall denote the partial action $\alpha|_{\G_{[e]}}$ of $\G_{[e]}$ on $\A_{[e]}$ by $\alpha_{[e]}$.
 	
 \end{notation}

 \subsection{Trace maps and subring of invariants}
 Let $\G$ be a finite connected groupoid with $\G_0=\{e_1,\ldots,e_r\}$, $\A=\oplus_{i=1}^r\A_i$ be a unital ring, $\A_i:=\A_{e_i}$ and $1_i$ be the identity element of $\A_i$ and $\alpha=(\A_g,\alpha_g)_{g\in \G}$ a unital partial action of $\G$ on $\A$.  \smallbreak

 Given $e_i,e_j\in \G_0$, we
 define the {\it trace map from $e_i$ to $e_j$} by:
 \begin{align*}
 t_{i,j}\colon\A\to \A_j,\quad t_{i,j}(a)=\sum_{g\in
 	\G(e_i,e_j)}\alpha_{g}(a1_{g^{-1}}),\quad a\in \A.
 \end{align*}
 Since $\G$ is connected we have that $\G(e_i,e_j)\neq \emptyset$. Moreover, $\A_g\unlhd \A_{t(g)}=\A_{j}$, for all $g\in \G(e_i,e_j)$. Hence, $t_{i,j}$ is well-defined. 
 
 Given $a\in \A$, consider the unique elements $a_k\in \A_k$ such that $a=\sum_{k=1}^ra_k$. Notice that 
 \[t_{i,j}(a)=\sum_{g\in
 	\G(e_i,e_j)}\alpha_{g}(a1_{g^{-1}})=\sum_{g\in 	\G(e_i,e_j)}\sum_{k=1}^r\alpha_{g}(a_k1_{g^{-1}}).\] Since $a_k1_{g^{-1}}\in \A_k\cap \A_i$, it follows that
 \begin{align}\label{trace_equal}
 t_{i,j}(a)=\sum_{g\in 	\G(e_i,e_j)}\alpha_{g}(a_i1_{g^{-1}})=t_{i,j}(a_i).
 \end{align}
 Note that $t_{i,i}$ 
 coincides with the partial trace map $t_{\alpha_{(e_i)}}$ introduced in   \cite[\S\,$2$]{DFP} for partial actions of groups.
 \smallbreak
 

 We also define the {\em subring of invariants of $\A$ under $\af$ from $e_i$ to $e_j$} by:
 \begin{align*}
 \A^{(i,j)}:=\{a\in \A\, : \,\alpha_{g}(a1_{g^{-1}})= a1_{g},\,\,\text{for all}\,\, g\in \G(e_i,e_j)\}.
 \end{align*}
 It is easy to check that $\A^{{(i,j)}}$ is a subring of $\A$. Also, $\A^{{(i,i)}}\cap \A_{i}=\A_i^{\af_{(e_i)}}$, where $\A_i^{\af_{(e_i)}}$ is the subring of invariants of $\A_i$ under the partial action $\alpha_{(e_i)}$ associated to the isotropy group $\G(e_i,e_i)$; see $\S\,$2 of \cite{DFP} for more details.
 
 \begin{prop1}\label{fixo} Let $e_i,e_j\in \p$, $x\in \A^{(i,j)}$ and $a\in \A$. Then the following statements hold:
 	\begin{enumerate}\renewcommand{\theenumi}{\roman{enumi}}   \renewcommand{\labelenumi}{(\theenumi)}
 		\item	$t_{i,j}(\A)\subseteq	\A_j^{\alpha_{(e_j)}}$;
 		\item  $t_{i,j}(xa)=xt_{i,j}(a)$ and $t_{i,j}(ax)=t_{i,j}(a)x$.
 		
 	\end{enumerate}	
 \end{prop1}
 
 \pf  Let $a\in \A$ and $h\in \G(e_j)$. Then, 
 \begin{align*}
 \alpha_h(t_{i,j}(a)1_{h^{-1}})=&\sum_{g\in \G(e_i,e_j)}\alpha_{h}(\alpha_{g}(a1_{g^{-1}})1_{h^{-1}})&\\
 \stackrel{\eqref{iguali}} {=}&\sum_{g\in \G(e_i,e_j)}\alpha_{hg}(a1_{(hg)^{-1}})1_h&\\
 =&\sum_{l\in \G(e_i,e_j)}\alpha_l(a1_{l^{-1}})1_h=t_{i,j}(a)1_h.&\
 \end{align*}
 Thus, $t_{i,j}(\A)\subseteq \A_j^{\alpha_{(e_j)}}$. The item \rm{(ii)} is immediate. 
 \epf
 
 Consider the maps $t_j:\A\to \A_j^{\alpha_{(e_j)}}$ given by
 \[t_j(a)=\sum\limits_{i=1}^r t_{i,j}(a), \text{ for all } a\in \A \text{ and } j=1,\ldots,r.\]
 By Proposition \ref{fixo} (i), $t_j$ is well-defined, for all $1\leq j\leq r$.
 
 Now we explore some properties of the maps  $t_j.$  Firstly, we recall that the trace map $t_{\alpha}:\A\to \A^{\alpha}$ associated to the partial action $\af$ was defined in \cite[p. 3668]{BP} by
 \[t_{\alpha}(a)=\sum_{g\in \G}\alpha_g(a1_{g^{-1}}),\quad a\in\A.\] 
 \begin{lem1} \label{talfa} Let $a\in \A$ and $g\in \G(e_i,e_j)$. Then:
 	\begin{enumerate}\renewcommand{\theenumi}{\roman{enumi}}\renewcommand{\labelenumi}{(\theenumi)}
 		\item $t_{\alpha}(a)=\sum_{j=1}^r t_{j}(a)$;\vone
 		\item $\alpha_g(t_i(a)1_{g^{-1}})=t_j(a)1_g$.\vone
 	\end{enumerate}
 \end{lem1}
 \pf Note that
 \[t_{\alpha}(a)=\sum\limits_{i,j=1}^r\sum\limits_{g\in \G(e_i,e_j)}\alpha_g(a1_{g^{-1}})=
 \sum\limits_{i,j=1}^r t_{i,j}(a)=\sum\limits_{j=1}^r t_j(a),\]
 and  \rm{(i)} follows. For \rm{(ii)}, let  $i,j,k\in\{1,\ldots,r\}$ and $g\in \G(e_i,e_j)$. Then
 \begin{align*}
 \alpha_g(t_{k,i}(a)1_{g^{-1}})=&\sum\limits_{h\in \G(e_k,e_i)}\alpha_g(\alpha_h(a1_{h^{-1}})1_{g^{-1}})=\sum\limits_{h\in \G(e_k,e_i)}\alpha_{gh}(a1_{(gh)^{-1}})1_g\\=&\sum\limits_{l\in \G(e_k,e_j)}\alpha_l(a1_{l^{-1}})1_g= \,\,t_{k,j}(a)1_g.
 \end{align*}
 Hence
 \begin{align*}
 \alpha_g(t_i(a)1_{g^{-1}})&=\sum\limits_{k=1}^r\alpha_g(t_{k,i}(a)1_{g^{-1}})=\sum\limits_{k=1}^rt_{k,j}(a)1_g=t_j(a)1_g.
 \end{align*}
 \epf

 \section{The Partial skew groupoid ring} 
 Let $\af=(\A_g,\af_g)_{g\in \G}$ be a partial action of a groupoid $\G$ on a ring $\A$. The {\it partial skew groupoid ring} was defined in \cite{BFP} as the direct sum $$\A\star_\af \G:=\bigoplus_{g\in \G}\A_g\de_g,$$
 where the $\de_g '$s are symbols, with the usual addition and the multiplication given by
 \[
 (a_g\de_g)(b_h\de_h)=\left\{\begin{array}{lc}
 \af_g(\af_{g\m}(a_g)b_h)\de_{gh}, & \text{if } (g,h)\in\G^2, \\
 0, & \text{otherwise},
 \end{array}\right.\]
 for all $g,h\in \G$, $a_g\in \A_g$ and $b_h\in\A_h$. When $\af$ is a global action, $\A\star_{\af}\G$ is called the {\it skew groupoid ring}.\vone
 
 It was observed in \cite{BFP} that this multiplication is well defined and  in general
 $\A\star_\af \G$ is neither associative nor unital. However, when $\af$ is unital, it is immediate to check that
 $\A\star_\af \G$ is associative. By    \cite[Proposition 3.3]{BFP}, if $\af$ is unital and $\p$ is finite, then $\A\star_\af \G$ is unital with $1_{\A\star_{\af} \G}=\sum_{e\in \p}1_e\delta_e$. 
 
 \begin{remark}
 	Let $\G$ be a groupoid such that $\p$ is finite, $\A$ a unital ring and $\af$ a unital partial action of $\G$ on $\A$. The map  $\varphi\colon \A\to \A\star_{\alpha}\G$, given by
 	\begin{align*}
 	\varphi(a)=\sum_{e\in \p}(a1_e)\delta_e, \quad a\in \A,
 	\end{align*}
 	is a ring monomorphism and $\varphi(1_\A)=1_{\A\star_{\alpha}\G}$. Consequently,  $\A\star_{\alpha}G$ is an $(\A,\A)$-bimodule  via:
 	\begin{align}\label{aabim}
 	&a\cdot (a_g\delta_g)=aa_g\delta_g,\quad\quad
 	(a_g\delta_g)\cdot a=a_g\alpha_g(a1_{g^{-1}})\delta_g,&\
 	\end{align}
 	for all $g\in\G$, $a_g\in \A_g$ and $a\in \A$. 
 \end{remark}
 
 \begin{teo1}\label{skew_component}
 	Let $\G$ be a  groupoid such that  $\G_0$ is finite, $\A$ a unital ring and $\alpha=(\A_g,\alpha_g)_{g\in \G}$ a unital partial action of $\G$ on $\A$. Assume that 
 	$\A=\oplus_{e\in \p}\A_e$. Then the following statements hold:
 	\begin{enumerate}\renewcommand{\theenumi}{\roman{enumi}}   \renewcommand{\labelenumi}{(\theenumi)}
 		\item $\A\star_\af \G=\bigoplus\limits_{[e]\in \p/\sim} \A_\ee\star_{\af_{[e]}}\G_{[e]}$  and $\B_{[e]}:=\A_\ee\star_{\af_{[e]}}\G_{[e]}$ is a two-sided ideal of $\A\star_\af \G$, for all $[e]\in \p/\hspace*{-0.15cm}\sim$. \vone
 		\item Let $u_i:=\sum_{f\in [e_i]}1_f\delta_f$ the identity element of $\B_{[e_i]}$. Then  $\{u_1,\ldots,u_n\}$ is a set of central pairwise orthogonal idempotents of $\A\star_\af \G$ such that 
 		\begin{align*}
 		1_{\A\star_\af \G}=\sum_{i=1}^n u_{i}\quad \text{and} \quad\B_{[e_i]}=(\A\star_\af \G) u_i.
 		\end{align*}
 		\item $\B_{[e]}\otimes_\A \B_{[f]}=0$, for all $[e],\,[f]\in \p/\hspace*{-0.15cm}\sim$ such that $[e]\neq [f]$.  \vone
 		
 		\item $\B_{[e]}\otimes_\A \B_{[e]}\simeq \B_{[e]}\otimes_{\A_{[e]}} \B_{[e]}$ as $(\A,\A)$-bimodules. \vone
 		
 		\item  $(\A\star_\af \G)\otimes_\A (\A\star_\af \G)\,\,\simeq\bigoplus\limits_{[e]\in \p/\sim}\B_{[e]}\otimes_{\A_\ee}\B_{[e]}$  as $(\A,\A)$-bimodules. 
 	\end{enumerate}
 \end{teo1}
 \pf \rm{(i)} By Proposition \ref{partial_action_comp_connected_2} (i), we can consider the partial skew groupoid ring $\A_\ee\star_{\af_{[e]}}\G_{[e]}$ which is a subring of $\A\star_\af \G$. Let $[e]\in \p/\hspace*{-0.15cm}\sim$, $g\in \G_{[e]}$, $h\in \G$, $a\in \A_g$ and $b\in \A_h$. If $(g,h)\notin \G^2$ then $(a\delta_g)(b\delta_h)=0$. If $(g,h)\in \G^2$ then $(a\delta_g)(b\delta_h)=a\af_g(b1_{g^{-1}})\delta_{gh}$. Observe that $a\af_g(b1_{g^{-1}})\in \A_{[e]}$, $t(gh)=t(g)\in [e]$ and $s(gh)=s(h)\sim t(h)=s(g)\in [e]$. Hence, $s(gh),t(gh)\in [e]$. Consequently, $(a\delta_g)(b\delta_h)\in \A_\ee\star_{\af_{[e]}}\G_{[e]}$. Similarly, we can verify that $\A_\ee\star_{\af_{[e]}}\G_{[e]}$ is a left ideal of $\A\star_\af \G$. Note also that 
 \begin{align*}
 \A\star_\af \G =\bigoplus_{g\in \G}\A_g\de_g=\bigoplus\limits_{[e]\in \p/\sim\,\,}\bigoplus_{g\in \G_{[e]}}\A_g\de_g=\bigoplus\limits_{[e]\in \p/\sim}\,\A_\ee\star_{\af_{[e]}}\G_{[e]}.\
 \end{align*}

 \noindent \rm{(ii)} It is clear that  $u_i$ is the identity element of $\B_{[e_i]}$. Using (i), we obtain that $1_{\A\star_\af \G}=\sum_{i=1}^n u_{i}$. Clearly, $\{u_1,\ldots,u_n\}$ is a set of central pairwise orthogonal idempotents of $\A\star_\af \G$ and $\B_{[e_i]}=(\A\star_\af \G) u_i$.\vone
 
 \noindent \rm{(iii)} Let $g\in \G_{[e]}$ and $h\in \G_{[f]}$. Then
 \begin{align*}
 \A_g\de_g\otimes_\A\A_h\de_h&=\A_g\de_g\otimes_\A 1_h\A_h\de_h\\
 &=\A_g\af_g(1_h1_{g\m})\de_g\otimes_\A \A_h\de_h.\
 \end{align*}
 Observe that $1_h1_{g\m}\in \A_h\cap \A_{g^{-1}}\subseteq \A_{[f]}\cap \A_{[e]}=\{0\}$. Thus, $\A_g\de_g\otimes_\A\A_h\de_h=0$. \vone
 
 \noindent \rm{(iv)} Let $g,h\in \G_{[e]}$. We define the following  map
 \begin{align*}
 &\varphi_{g,h}:\A_g\de_g\times\A_h\de_h\to \A_g\de_g\otimes_{\A_{[e]}}\A_h\de_h,& &\varphi_{g,h}(x\de_g,y\de_h)=x\de_g\otimes_{\A_{[e]}}y\de_h,& 
 \end{align*}
 for all $x\in \A_{g}$ and $y\in \A_{h}$. Given $a=\sum_{f\in\p}a_f\in \A$, note that 
 \begin{align*}
 \varphi_{g,h}(x\de_g\cdot a,y\de_h) &=\varphi_{g,h}(x\alpha_g(a1_{g^{-1}})\de_g,y\de_h)\\
 &=x\alpha_g(a_{s(g)}1_{g^{-1}})\de_g\otimes_{\A_{[e]}} y\de_h\\
 &=x\de_g\cdot a_{s(g)} \otimes_{\A_{[e]}} y\de_h\\
 &=x\de_g \otimes_{\A_{[e]}} a_{s(g)}\cdot  y\de_h\\
 &=x\de_g \otimes_{\A_{[e]}} a_{t(h)}y\de_h\\
 &=x\de_g \otimes_{\A_{[e]}} ay\de_h\\
 &=\varphi_{g,h}(x\de_g, a\cdot y\de_h).\
 \end{align*}
 Consequently, 
 \begin{align*}
 &\overline{\varphi}_{g,h}:\A_g\de_g\otimes_\A\A_h\de_h\to \A_g\de_g\otimes_{\A_{[e]}}\A_h,& &\overline{\varphi}_{g,h}(x\de_g\otimes_\A y\de_h)=x\de_g\otimes_{\A_{[e]}}y\de_h,& 
 \end{align*}
 is well-defined map and it is clear that  $\overline{\varphi}_{g,h}$ is an $(\A,\A)$-bimodule homomorphism. 
 Similarly, we can prove that 
 \begin{align*}
 \overline{\psi}_{g,h}:\A_g\de_g\otimes_{\A_{[e]}}\A_h\de_h\to \A_g\de_g\otimes_{\A}\A_h,& &\overline{\psi}_{g,h}(x\de_g\otimes_{\A_{[e]}} y\de_h)=x\de_g\otimes_{\A}y\de_h,& \
 \end{align*}
 is a well-defined $(\A,\A)$-bimodule homomorphism. Note that $\overline{\psi}_{g,h}$ is the inverse of $\overline{\varphi}_{g,h}$, for all $(g,h)\in \G^2$. 
 The family of isomorphisms $\{\overline{\varphi}_{g,h}\}_{g,h\in \G_{[e]}}$ induces the isomorphism required.
 \vone
 
 \noindent \rm{(v)} Follows from \rm{(i)}, \rm{(iii)} and \rm{(iv)}.   \epf

 \section{Separability of the partial skew groupoid ring}
 
 Let $\G$ be a groupoid such that $\G_0$ is finite, $\A$ a unital ring and $\alpha=(\A_g,\alpha_g)_{g\in \G}$ a unital partial action of $\G$ on $\A$. We are interested in determining when the ring extension $\A\subset \A\star_\af \G$ is separable. This problem  can be reduced to the case where $\G$ is a connected groupoid, as we shall see in the next result.

 \begin{prop1}\label{prop:redu} Let $\G$ be a  groupoid such that $\G_0$ is finite and $\A$ be a unital ring. Suppose that $\alpha=(\A_g,\alpha_g)_{g\in \G}$ is a unital partial action of $\G$ on $\A$ and 
 	$\A=\oplus_{e\in \p}\A_e$. Then $\A\subset \A\star_\af \G $ is separable if and only if $\A_{[e]}\subset \A_{[e]}\star_{\af_{[e]}} \G_{[e]}$ is separable, for all  $[e]\in \p/\hspace*{-0.15cm}\sim$.
 \end{prop1}
 \pf Assume that $\B:=\A\star_\af \G\supseteq \A$ is a separable extension and take $x\in \B\otimes_\A \B$ an idempotent of separability. Suppose that $x=\sum_{k=1}^r b_k\otimes_{\A}b'_k$, with $b_k,\,b'_k\in \B$. By Theorem \ref{skew_component} (i), there exist elements $b_{k,i},\,b'_{k,i}\in \B_{[e_i]}$ such that $b_k=\sum_{i=1}^n b_{k,i}$ and $b'_k=\sum_{i=1}^n b'_{k,i}$. Fix $x_i:=\sum_{k=1}^{r}b_{k,i}\otimes_{\A_{[e_i]}} b'_{k,i}\in \B_{[e_i]}\otimes_{\A_{[e_i]}} \B_{[e_i]}$, for all $1\leq i\leq n$. Note that
 \begin{align*}
 1_{\B}&=m_{\B}(x)=\sum\limits_{k=1}^{r}\sum\limits_{i,j=1}^{n}b_{k,i}b'_{k,j}=\sum\limits_{k=1}^{r}\sum\limits_{i=1}^{n}b_{k,i}b'_{k,i}\\
 &=\sum\limits_{i=1}^{n}\sum\limits_{k=1}^{r}b_{k,i}b'_{k,i}=\sum\limits_{i=1}^{n}m_{\B_{[e_i]}}(x_i).
 \end{align*}
 By Theorem \ref{skew_component} (ii),  $m_{\B_{[e_i]}}(x_i)=u_i=1_{\B_{[e_i]}}$. Given an element $b\in \B_{[e_i]}$, we have  $bx=\sum_{k=1}^rbb_{k,i}\otimes_{\A}b'_{k,i}=\sum_{k=1}^rb_{k,i}\otimes_{\A}b'_{k,i}b=xb$. The isomorphism given in Theorem \ref{skew_component} (vi) implies    
 $bx_i=\sum_{k=1}^rbb_{k,i}\otimes_{\A_{[e_i]}}b'_{k,i}=\sum_{k=1}^rb_{k,i}\otimes_{\A_{[e_i]}}b'_{k,i}b=x_ib$. Therefore, $x_i$ is an idempotent of separability of $\B_{[e_i]}$, for all $i=1,\ldots,n$. 
 
 Conversely, suppose that $\B_{[e_i]}=\A_{[e_i]}\star_{\af_{[e_i]}} \G_{[e_i]}\supseteq\A_{[e_i]}$ is a separable extension, for each $i=1,\ldots,n$. Consider $x_{i}\in \B_{[e_i]}\otimes_{\A_{[e_i]}} \B_{[e_i]}$ an idempotent of separability of $\B_{[e_i]}$. Suppose that $x_i=\sum_{k=1}^{r_i}b_{k,i}\otimes_{\A_{[e_i]}} b'_{k,i}$, where $b_{k,i}$ and $b'_{k,i}$ are elements of $\B_{[e_i]}$, for all $k=1,\ldots,r_i$. It is straightforward to verify that 
 $x:=\sum_{i=1}^{n}\sum_{k=1}^{r_i}b_{k,i}\otimes_{\A}b'_{k,i}$ is an idempotent of separability of $\A\star_\af \G$.\epf

 In what follows in this section, $\G$ is a connected groupoid with  finite set of objects $\G_0=\{e_1,\ldots,e_r\}$, $\A=\oplus_{i=1}^r\A_i$ is a unital ring where $\A_i=\A_{e_i}$, $\alpha=(\A_g,\alpha_g)_{g\in \G}$ is a unital partial action of $\G$ on $\A$ and  $1_i$ is the identity element of $\A_i$.

 \begin{teo1}\label{sepa} If $\G$ is finite then the following statements are
 	equivalent:
 	\begin{enumerate}\renewcommand{\theenumi}{\roman{enumi}}   \renewcommand{\labelenumi}{(\theenumi)}
 		\item $\A\subset \A\star_{\alpha}\G$ is a separable extension;\vone
 		\item there exists $a\in C(\A)$ such that $t_i(a)=1_i$, for all $i=1,\ldots, r$.
 	\end{enumerate}
 \end{teo1}
 
 \pf Assume that $\A\subset \B:=\A\star_{\alpha}\G$ is a  separable
 extension of rings. Consider an idempotent of separability $x\in \B\otimes_\A \B$. Given $(g,h)\in \G^2$, $a\in \A_g$ and $b\in \A_h$ we have
 $(a\delta_g)\otimes_{\A} (b\delta_h)=(a\delta_g)\cdot b\otimes_{\A} 1_h\delta_h=a\alpha_g(b1_{g^{-1}})\delta_g\otimes_{\A}
 1_h\delta_h$. Thus we can suppose that \[x=\sum\limits_{(g,h)\in
 	\G^2}a_{g,h}\delta_g\otimes_{\A} 1_h\delta_h,\] with $a_{g,h} \in \A_g\cap \A_{gh}$.
 Since $m_{\B}(x)=1_{\B}$, we have
 \[\sum\limits_{(g,h)\in \G^2}a_{g,h}\delta_{gh} =\sum\limits_{i=1}^r1_i\delta_{e_i}.\] 
 Let $\G^{2}(l):=\{(g,h)\in \G^2\,:\,gh=l\}$, for all $l\in \G$. Then
 \begin{align}\label{suma}
 \sum\limits_{(g,h)\in \G^2(l)}a_{g,h}=\left\{\begin{array}{l}
 1_i,\quad\text{if}\quad l=e_i; \\
 0,\quad\,\,\text{if}\quad l \notin \G_0. \\
 \end{array}\right.
 \end{align}
 Since $x$ is $\B$-central, we get
 \begin{align*}
 \sum\limits_{(g,h)\in \G^2}aa_{g,h}\delta_g\otimes_{\A} 1_h\delta_h
 &=\sum\limits_{(g,h)\in \G^2}a_{g,h}\delta_g\otimes_{\A} (1_h\delta_h)\cdot a\\
 &\stackrel{\eqref{aabim}}{=}\sum\limits_{(g,h)\in \G^2}a_{g,h}\delta_g\otimes_{\A}\alpha_h(a1_{h^{-1}})\delta_h\\
 & =\sum\limits_{(g,h)\in \G^2}a_{g,h}\alpha_g(\alpha_h(a1_{h^{-1}})1_{g^{-1}})\delta_g\otimes_{\A} 1_h\delta_h,\
 \end{align*}
 for all $a\in\A$. Particularly, if we consider
 $g=h=e_i$  then $aa_{e_i,e_i}=a_{e_i,e_i}a$, for all
 $a\in \A$. Thus, $a_{e_i,e_i}\in C(\A)$, for all $i=1,\ldots,r$. Observe
 also that
 \begin{align}\label{lx}
 (1_l\delta_l)x=\sum\limits_{{(l,g),(g,h)\in
 		\G^2}}\alpha_{l}(a_{g,h}1_{l^{-1}})\delta_{lg}\otimes_{\A} 1_h\delta_h,
 \end{align}
 and
 \begin{align}\label{xl}
 x(1_l\delta_l)=\sum\limits_{{(g,h),(h,l)\in \G^2}}a_{g,h}\delta_g\otimes_{\A} 1_h1_{hl}\delta_{hl},
 \end{align}
 for all $l\in \G$. Since $(1_l\delta_l)x=x(1_l\delta_l)$, comparing
 the coordinates of $\delta_g\otimes_{\A} \delta_{s(g)}$ in (\ref{lx}) and (\ref{xl}) we obtain that $\alpha_g(a_{s(g),s(g)}1_{g^{-1}})\delta_g\otimes_{\A} 1_{s(g)}\delta_{s(g)}=a_{g,g^{-1}}\delta_g\otimes_{\A} 1_{s(g)}\delta_{s(g)},$ which implies
 \begin{align}\label{iguald}
 \alpha_g(a_{s(g),s(g)}1_{g^{-1}})=a_{g,g^{-1}}.
 \end{align}
 Take $a:=a_{e_1,e_1}+\ldots+a_{e_r,e_r}\in C(\A)$. Notice that
 \begin{align*}
 t_i(a)&=\sum_{k=1}^r t_{k,i}(a)\stackrel{\eqref{trace_equal}}{=}\sum_{k=1}^r t_{k,i}(a_{e_k,e_k})=\sum_{k=1}^r\sum\limits_{g\in \G(e_k,e_i)}\alpha_g(a_{e_k,e_k}1_{g^{-1}})\\
 &\stackrel{\eqref{iguald}}{=}\sum_{k=1}^r\sum\limits_{g\in \G(e_k,e_i)} a_{g,g^{-1}}=\sum_{g\in \G\atop t(g)=e_i}a_{g,g^{-1}}\stackrel{\eqref{suma}}{=}1_i,\
 \end{align*}
 for all $i=1,\ldots,r$.\vone
 
 Conversely, let $a\in C(\A)$ such that $t_i(a)=1_i$, for all $i=1,\ldots,r$. Suppose that $a=a_1+\ldots+a_r$, with $a_i\in \A_i$.
 Then, we take the following element
 \[x:=\sum\limits_{i,j=1}^r\sum\limits_{g\in \G(e_j,e_i)}\alpha_g(a_j1_{g^{-1}})\delta_g\otimes_{\A}
 1_{g^{-1}}\delta_{g^{-1}}\in \B\ot_\A \B.\] Note that
 \begin{align*}
 m_{\B}(x)&=\sum\limits_{i=1}^r\sum\limits_{j=1}^r\sum\limits_{g\in \G(e_j,e_i)}\alpha_g(a_j1_{g^{-1}})\delta_{e_i}\\
 &=\sum\limits_{i=1}^r\sum\limits_{j=1}^rt_{j,i}(a_j)\delta_{e_i}\\&\stackrel{\eqref{trace_equal}}{=}
 \sum\limits_{i=1}^r\sum\limits_{j=1}^rt_{j,i}(a)\delta_{e_i}\\
 &=\sum\limits_{i=1}^r t_i(a)\delta_{e_i}\\&=\sum\limits_{i=1}^r 1_i\delta_{e_i}=1_{\B}.
 \end{align*}
 Since $a\in C(\A)$, it follows that $a_i\in C(\A)$, for all $i=1,\ldots,n$. Consequently,
 $\alpha_g(a_i1_{g^{-1}})\in C(\A)$, for all $i=1,\ldots,n$. In order to prove that $x$ is $\B$-central it is enough to show that  $xb=bx$ and $x(1_l\delta_l)=(1_l\delta_l)x$, for 
 all $b\in \A$ and $l\in \G$. Given $b\in \A$  and $l\in \G$, we have that
 \begin{align*}
 xb&= \sum\limits_{i,j=1}^r\sum\limits_{g\in \G(e_i,e_j)}\alpha_g(a_i1_{g^{-1}})\delta_g\otimes_{\A}
 \af_{g^{-1}}(b1_g)\delta_{g^{-1}}\\
 &= \sum\limits_{i,j=1}^r\sum\limits_{g\in \G(e_i,e_j)}\alpha_g(a_i1_{g^{-1}})b\delta_g\otimes_{\A}
 1_{g^{-1}}\delta_{g^{-1}}\\
 &=\sum\limits_{i,j=1}^r\sum\limits_{g\in \G(e_i,e_j)}b\alpha_g(a_i1_{g^{-1}})\delta_g\otimes_{\A}
 1_{g^{-1}}\delta_{g^{-1}}\\
 &=bx,
 \end{align*}
 and
 \begin{align*}
 (1_l\delta_l)x&=\sum\limits_{i,j=1}^r\sum\limits_{g\in \G(e_i,e_j)}(1_l\delta_l)(\alpha_g(a_i1_{g^{-1}})\delta_g)\otimes_{\A}
 1_{g^{-1}}\delta_{g^{-1}}\\
 &=\sum\limits_{i=1}^r\sum\limits_{g\in \G(e_i,s(l))}
 \alpha_{lg}(a_i1_{(lg)^{-1}})1_l\delta_{lg}\otimes_{\A}
 1_{g^{-1}}\delta_{g^{-1}}\quad \quad\quad (lg=h)\\
 &=\sum\limits_{i=1}^r\sum\limits_{h\in \G(e_i,t(l))}
 \alpha_h(a_i1_{h^{-1}})\alpha_h(1_{h^{-1}}1_{h^{-1}l})\delta_h
 \otimes_{\A} 1_{h^{-1}l}\delta_{h^{-1}l}\\
 &=\sum\limits_{i=1}^r\sum\limits_{h\in \G(e_i,t(l))}
 \alpha_h(a_i1_{h^{-1}})\delta_h\cdot 1_{h^{-1}l}
 \otimes_{\A} 1_{h^{-1}l}\delta_{h^{-1}l}\\
 &=\sum\limits_{i=1}^r\sum\limits_{h\in \G(e_i,t(l))}
 \alpha_h(a_i1_{h^{-1}})\delta_h
 \otimes_{\A} 1_{h^{-1}l}\delta_{h^{-1}l}\\
 &=\sum\limits_{i=1}^r\sum\limits_{h\in \G(e_i,t(l))}
 \alpha_h(a_i1_{h^{-1}})\delta_h\otimes_{\A}
 (1_{h^{-1}}\delta_{h^{-1}})(1_l\delta_l)\\
 &=x(1_l\delta_l).\
 \end{align*}\epf

 Now we present two  examples.
 
 \begin{exe1}\label{psepa} Let $\G=\{g,g^{-1}\}$ be the groupoid with $e_1=s(g)\neq t(g)=e_2$.
 	Consider  the ring $\A=\oplus_{i=1}^4\K v_i$, where
 	$v_i$ are pairwise orthogonal idempotents with sum $1_\A$ and $\Bbbk$ is a field. We put
 	$\A_{s(g)}=\K v_1\oplus \K v_2$, $\A_{t(g)}=\K v_3\oplus \K
 	v_4$, $\A_{g^{-1}}=\K v_2$ and $\A_{g}=\K v_3$. Defining $\af_{g}(v_2)=v_3$, $\af_{g^{-1}}(v_3)=v_2$, $\af_{s(g)}={\rm
 		id}_{\A_{s(g)}}$ and $\af_{t(g)}={\rm
 		id}_{\A_{t(g)}}$, we obtain that $\alpha=(\A_g,\alpha_g)_{g\in \G}$ is a unital partial action of $\G$ on $\A$.
 	Given an element $a=\sum_{i=1}^4\lambda_iv_i\in \A$, we have that \[t_1(a)=\lambda_1v_1+(\lambda_2+\lambda_3)v_2, \quad t_2(a)=(\lambda_2+\lambda_3)v_3+\lambda_4v_4.\]
 	Consequently, $t_1(a)=v_1+v_2=1_1$ and $t_2(a)=v_3+v_4=1_2$ if and only if $a=v_1+\lambda v_2+(1-\lambda)v_3+v_4$, where $\lambda\in \K$. By Theorem \ref{sepa}, the extension  
 	$\A\subset \A\star_{\alpha}\G$ is separable. Moreover, fixing $a_1=v_1+\lambda v_2$ and $a_2=(1-\lambda)v_3+v_4$ we have that
 	\[x_{\lambda}:=a_1\delta_{s(g)}\otimes 1_1\delta_{s(g)}+\lambda v_3\delta_g \otimes v_2\delta_{g^{-1}} + (1-\lambda)v_2\delta_{g^{-1}} \otimes v_3\delta_{g}   +a_2\delta_{t(g)}\otimes 1_2\delta_{t(g)}\]
 	is an idempotent of separability of $\A\star_{\alpha}\G$, for all $\lambda\in \K$.
 \end{exe1}
 
 \begin{exe1} Consider the groupoid $\G=\{g,h,x,y,x^{-1},y^{-1}\}$ with set of objects $\G_0=\{e_1,e_2\}$ and satisfying 
 	\begin{align*}
 	s(x)=s(y)=e_1&,\quad  t(x)=t(y)=e_2,\quad s(g)=t(g)=e_1,\quad s(h)=t(h)=e_2,&\\
 	& g^2=e_1,\quad \quad h^2=e_2,\quad\quad xg=y=hx.&
 	\end{align*}	
 	Let $\A=\oplus_{i=1}^2\K v_i$ be the ring, where
 	$v_i$ are pairwise orthogonal idempotents with sum $1_\A$ and $\Bbbk$ is a field. We put
 	$\A_{e_1}=\A_{g}=\K v_1$, $\A_{e_2}=\A_{h}=\K v_2$, $\A_{x^{-1}}=\A_{x}=\{0\}$ and $\A_{y^{-1}}=\A_{y}=\{0\}$. 
 	Defining $\af_{x}=\af_{x^{-1}}=\af_{y}=\af_{y^{-1}}=0$, $\af_{e_1}=\af_{g}=\Id_{\A_{e_1}}$ and $\af_{e_2}=\af_{h}=\Id_{\A_{e_2}}$, it follows that $\alpha=(\A_l,\alpha_l)_{l\in \G}$ is a unital partial action of $\G$ on $\A$.
 	Given an element $a=\lambda_1v_1+\lambda_2v_2\in \A$, we have that $t_1(a)=2\lambda_1v_1$ and $t_2(a)=2\lambda_2v_2$. Thus, $\A\subset \A\star_{\alpha}\G$ is separable if and only if the characteristc of $\K$ is not $2$.
 \end{exe1}

\begin{remark}
 	Let $G$ be a finite group, $\B$ a unital ring and $\gamma=(\B_g,\gamma_g)_{g\in G}$ a unital partial action of $G$ on $\B$. By  \cite[Theorem 3.1 (ii)]{BLP}, 
 	$\B\subset \B\star_{\gamma}G$ is a separable extension if and only if there exists $b\in C(\B)$ such that $t_{\gamma}(b)=1_{\B}$. Since $t_{i}=t_{\gamma}$,  then Theorem \ref{sepa} generalizes this result to the context of partial groupoid actions. 
 \end{remark}
 
 \begin{cor1} \label{cor:sepa} Suppose that $\G$ is finite and that
 	\begin{itemize}
 		\item [$\circ$] there exist $a\in C(A)$ and $k\in\{1,\ldots,r\}$ such that $t_k(a)=1_k$,
 		\item [$\circ$] for each $j\in\{1,\ldots,r\}$, there exists $g_j\in \G(e_k,e_j)$ such that $1_{g_j}=1_j$.
 	\end{itemize}
 	Then $\A\subset \A\star_{\alpha}\G$ is a separable extension.
 \end{cor1}
 \pf By Lemma \ref{talfa} (ii), 
 \[1_{g_j}=\af_{g_j}(1_k1_{g_j^{-1}})=\af_{g_j}(t_k(a)1_{g_j^{-1}})=t_j(a)1_{g_j}.\]
 Thus $0=(t_j(a)-1_j)1_{g_j}=(t_j(a)-1_j)1_{j}=t_j(a)-1_j$ and whence the result follows from Theorem \ref{sepa}.
 \epf
 
 We recall from \cite{HS} that an extension of rings $\A\subset \B$ is said  semisimple if any
 left $\B$-submodule $W$ of a left $\B$-module $V$ having an $\A$-complement
 in $V$, has a $\B$-complement in $V$. Equivalently, $\A\subset \B$ is semisimple if any exact sequence in the category $_\B\mathcal{M}$ of left $\B$-modules that splits in $_\A\mathcal{M}$ also splits in $_\B\mathcal{M}$. It is well-known that separable ring extensions are semisimple; see e.~g.   \cite[Remark 1.1]{CIM}. Thus we have the following immediate consequence of Theorem \ref{sepa}.

 \begin{cor1} If $\G$ is finite and there exists $a\in C(\A)$ such that $t_i(a)=1_i$, for all $i=1,\ldots, r$, then 
 	$\A\subset \A\star_{\alpha}\G$ is a ring semisimple extension. \qed
 \end{cor1}
\section{Separability of the  skew groupoid ring}

Throughout this section, $\G$ is a connected groupoid such that $\G_0=\{e_1,\ldots,e_r\}$, $\A=\oplus_{i=1}^r\A_i$ is a unital ring where $\A_i=\A_{e_i}\neq 0$, $\alpha=(\A_g,\alpha_g)_{g\in \G}$ is a unital global action of $\G$ on $\A$ and  $1_i$ is the identity element of $\A_i$.

If $\af$ is a partial action (not global) then the separability of $\A\star_{\alpha}\G$ over $\A$ does not implies the finiteness of $\G$, even when $\G$ is a group; see Remark 3.2 of \cite{BLP}. However,  the statement is true for global actions as we shall see below.  

\begin{prop1}\label{csepa} If $\A\subset \A\star_{\alpha}\G$ is a
	separable extension then $\G$ is finite.
\end{prop1}

\pf Fix $\B:=\A\star_{\alpha}\G$. Let $x\in \B\otimes_{\A} \B$ an idempotent of separability. As in the proof of Theorem \ref{sepa}, we can take $x=\sum_{(g,h)\in \G^2}a_{g,h}\delta_g\otimes 1_h\delta_h$, where $a_{g,h}\in \A_g\cap \A_{gh}$. Since $\af$  is global,
$\A_{g^{-1}}=\A_{s(g)}$ and whence $1_{g^{-1}}=1_{s(g)}$, for all $g\in \G$.
As in the proof of Theorem \ref{sepa}, it follows from (\ref{suma}) that 
\begin{align*}\label{1e}
1_i=\sum\limits_{g\in
	\G(e_j,e_i)}\alpha_g(a_{e_j,e_j}), \quad \text{ for all } 1\leq i,j\leq r.
\end{align*}
Since $\A_i\neq 0$, we obtain that $a_{e_j,e_j}\neq 0$. Therefore $a_{g,g^{-1}}=\alpha_{g}(a_{e_j,e_j})\neq 0$, for all
$g\in \G(e_j,e_i)$. Hence $\G(e_j,e_i)$ is finite, for all $1\leq i,j\leq n$, and whence $\G$ is finite. 
\epf

\begin{cor1}\label{cor:globalcase} The following statements are
	equivalent:
	\begin{enumerate}\renewcommand{\theenumi}{\roman{enumi}}   \renewcommand{\labelenumi}{(\theenumi)}
		\item $\A\subset \A\star_{\alpha}\G$ is a separable extension;\smallbreak
		\item $\G$ is finite and  there exists $a\in C(\A)$ such that $t_i(a)=1_i$, for all $i=1,\ldots, r$.
	\end{enumerate}
\end{cor1} 
\pf It follows from Proposition \ref{csepa} and Theorem \ref{sepa}. \epf

\begin{cor1}\label{cor:maisum}
	The following statements are
	equivalent:
	\begin{enumerate}\renewcommand{\theenumi}{\roman{enumi}}   \renewcommand{\labelenumi}{(\theenumi)}
		\item $\A\subset \A\star_{\alpha}\G$ is a separable extension;\smallbreak
		\item $\G$ is finite and  there exist $a\in C(\A)$  and $k\in\{1,\ldots,r\}$ such that $t_k(a)=1_k$.
	\end{enumerate}
\end{cor1} 
\pf \rm{(i)} $\Rightarrow $ \rm{(ii)} is immediate from Corollary \ref{cor:globalcase}. For the converse, notice that for each $g\in \G(e_k,e_j)$ we have that $1_j=1_{t(g)}=1_g$.  Hence, the result follows directly from Corollary \ref{cor:sepa}.
\epf

\begin{teo1}\label{teo:global}
	If $\A\subset \A\star_{\alpha}\G$ is a separable extension then there exist $1\leq i\leq r$ such that $\A_i\subset \A_i\star_{\alpha_{(e_i)}}\G(e_i)$ is a separable extension. 
\end{teo1}

\pf Since $\G$ is connected we can fix $g_k\in \G(e_i,e_k)$, for all $1\leq k\leq r$. Assume that $\A\subset \A\star_{\alpha}\G$ is separable. By Corollary \ref{cor:maisum}, $\G$ is finite and there exist an element  $b=\sum_{k=1}^{r}b_k\in C(\A)$, where $b_k\in C(\A_k)$, and $1\leq i\leq r$ such that $t_i(b)=1_i$. Taking $a=\sum_{k=1}^r\alpha_{g_k^{-1}}(b_k)\in C(\A_i)$, we obtain
\begin{align*}
t_{i,i}(a)=\sum\limits_{k=1}^{r}t_{i,i}(\alpha_{g_k^{-1}}(b_k))=\sum\limits_{k=1}^{r}t_{k,i}(b_k)=\sum\limits_{k=1}^{r}t_{k,i}(b)=t_i(b)=1_i.
\end{align*}
Hence, Theorem 3.1 (i) of \cite{BLP} implies that $\A_i\star_{\alpha_{(e_i)}}\G(e_i)$ is separable. \vone \epf

\begin{remark}
	Let  $e_i, e_j\in \G_0$ and fix $l\in \G(e_i,e_j)$. Given $a\in \A_i$ and $g\in \G(e_i)$, we define $\psi:\A_i\star_{\af_{(e_i)}}\G(e_i)\to\A_j\star_{\af_{(e_j)}}\G(e_j)$ by
	$\psi(a\delta_g)=\af_{l}(a)\delta_{lgl^{-1}}$. Extending by linearity we obtain that $\psi$ is a ring isomorphism. Thus, $\A_j\subset \A_j\star_{\alpha_{(e_j)}}\G(e_j)$ is separable for all $1\leq j\leq r$ if and only if there exist $1\leq i\leq r$ such that $\A_i\subset \A_i\star_{\alpha_{(e_i)}}\G(e_i)$ is separable.
\end{remark}


\begin{thebibliography}{9999}
 	
 	\bibitem{AS} R. Alfaro and G. Szeto, {\it The Centralizer on H-Separable Skew Group Rings.}  Rings, 
 	Extension and Cohomology, Marcel Dekker. \textbf{159}, (1994).
 	
 	\bibitem{AS1} R. Alfaro and G. Szeto, {\it Skew group rings which are Azumaya.}
 	Comm. Alg. \textbf{23}, 2255--2261 (1995).
 	
 	
 	
 	\bibitem{BFP} D. Bagio, D. Fl\^ores and A. Paques, {\it Partial actions of ordered groupoids on rings.}   J. Algebra Appl. \textbf{9},
 	501--517 (2010).
 	
 	\bibitem{BLP} D. Bagio, J. R. Lazzarin and A. Paques, {\it Crossed products by twisted partial actions: separability,
 		semisimplicity and Frobenius properties.}  Comm.  Alg.
 	\textbf{38}, 496--508 (2010).
 	
 	
 	\bibitem{BP} D. Bagio and A. Paques, {\it Partial Groupoid Actions: Globalization, Morita Theory, and Galois Theory.}
 	Comm. Alg. \textbf{40}, 3658--3678 (2012).
 	
 	\bibitem{CAS}
 	F. Casta\~{n}o, J. G\'{o}mez Torrecillas and C. N{\v a}st{\v a}sescu, {\it Separable functors in graded rings},
 	J. Pure Appl. Algebra {\bf 127}, 219--230 (1998).
 	
 	
 	\bibitem{CIM}Caenepeel, S., Ion, B., Militaru, G. {\it The structure of Frobenius algebras and
 		separable algebras.}  K-Theory \textbf{ 19} 365–402 (2000).
 	
 	\bibitem{DEM} F. M. DeMeyer and E. Ingraham, {\it Separable Algebras Over Commutative Rings}, Springer Lecture Notes 181 (1971).
 	
 	\bibitem{DFP} M. Dokuchaev, M. Ferrero and A. Paques, {\it Partial actions and Galois theory.}   J. Pure Appl. Algebra \textbf{208},
 	77--87 (2007).
 	
 	\bibitem{HS} K. Hirata and K. Sugano, {\it On semisimple extensions and separable extensions over noncommutative rings.}  J. Math Soc. Japan \textbf{18},
 	360--373 (1966).
 	
 	\bibitem{I} S. Ikehata, {\it Note on Azumaya Algebras and H-Separable Extensions.} Math. J. Okayama Univ. \textbf{23}, 17--18 (1981).
 	
 	
 	\bibitem{L1} P. Lundstr\"{o}m, {\it Crossed Product Algebras Defined by Separable Extensions.}  J. Algebra {\bf 283}, 723--737 (2005).
 	
 	\bibitem{L2} P. Lundstr\"om, {\it Separable groupoid rings.}  Comm. Algebra \textbf{34} (9), 3029--3041 (2005).

 	\bibitem{NOP} P. Nystedt, J. \"{O}inert, H. Pinedo, {\it Epsilon-strongly Graded Rings, Separability and Semisimplicity.} Preprint.

 	
 \end{thebibliography}
\end{document}